\documentclass[aps,superscriptaddress,preprintnumbers, superscriptaddress, nofootinbibt,onecolumn]{revtex4}
\usepackage{amsmath}
\usepackage{eurosym}
\usepackage{amssymb}
\usepackage{graphicx}
\usepackage{color}

\def\be{\begin{equation}}
\def\ee{\end{equation}}
\def\bea{\begin{eqnarray}}
\def\eea{\end{eqnarray}}

\begin{document}

\title{On the integrability of the Abel and of the extended Li\'{e}nard
equations}
\author{Man  Kwong Mak}
\email{mankwongmak@gmail.com}
\affiliation{Departamento de F\'{\i}sica, Facultad de Ciencias Naturales, Universidad de
Atacama, Copayapu 485, Copiap\'o, Chile}
\author{Tiberiu Harko}
\email{tiberiu.harko@gmail.com}
\affiliation{Department of Physics, Babes-Bolyai University, Kogalniceanu Street,
Cluj-Napoca 400084, Romania,}
\affiliation{School of Physics, Sun Yat-Sen University, Xingang Road, Guangzhou 510275, P. R. China}

\date{\today }

\begin{abstract}
We present some exact integrability cases of the extended Li\'{e}nard
equation $y^{\prime \prime }+f\left( y\right) \left(y^{\prime
}\right)^{n}+k\left( y\right) \left(y^{\prime }\right)^{m}+g\left(y\right)
y^{\prime }+h\left( y\right) =0$, with $n>0$ and $m>0$ arbitrary constants, while $f(y)$, $k(y)$, $g(y)$, and $h(y)$
are arbitrary functions. The solutions are obtained by transforming the equation Li\'{e}nard equation to an
equivalent first kind first order Abel type equation given by $\frac{dv}{dy}%
=f\left( y\right) v^{3-n}+k\left( y\right) v^{3-m}+g\left( y\right)
v^{2}+h\left( y\right) v^{3}$, with $v=1/y^{\prime }$. As a first step in
our study we obtain three integrability cases of the extended
quadratic-cubic Li\'{e}nard equation, corresponding to $n=2$ and $m=3$, by
assuming that particular solutions of the associated Abel equation are
known. Under this assumption the general solutions of the Abel and Li\'{e}%
nard equations with coefficients satisfying some differential conditions can
be obtained in an exact closed form. With the use of the Chiellini
integrability condition, we show that if a particular solution of the Abel
equation is known, the general solution of the extended quadratic cubic Li%
\'{e}nard equation can be obtained by quadratures. The Chiellini
integrability condition is extended to generalized Abel equations with $%
g(y)\equiv 0$ and $h(y)\equiv 0$, and arbitrary $n$ and $m$, thus allowing
to obtain the general solution of the corresponding Li\'{e}nard equation.
The application of the generalized Chiellini condition to the case of the
reduced Riccati equation is also considered.

{\bf Keywords:} Abel equation; Li\'{e}nard equation; general solutions; integrability conditions

{\bf MSC:}  34A05;  34A25; 34B30; 34C15; 34G20

  \end{abstract}

\maketitle

\tableofcontents


\section{Introduction}

The Li\'{e}nard type ordinary second order nonlinear differential equation
of the form \cite{Lien,Lien1}
\begin{equation}
\ddot{x}(t)+f(x)\dot{x}(t)+g(x)=0,
\end{equation}
as well as its generalization, the Levinson-Smith type equation \cite{Lev}
\begin{equation}  \label{Lev}
\ddot{x}(t)+f\left(x,\dot{x}\right)\dot{x}(t)+g(x)=0,
\end{equation}
where a dot represents the derivative with respect to the independent
variable $t$, $g$ is an arbitrary $C^1$ functions of $x$, while $f$ is an
arbitrary function of $x$ and $\dot{x}$, play an important role in many
areas of physics, biology and engineering \cite{2}. These types of equations
have been intensively studied from both mathematical and physical point of
view, and these investigation still remain an active and quickly developing
field of research in pure mathematics and mathematical physics \cite{lit0} -
\cite{lit6}.

One of the important applications of the generalized Li\'{e}nard type
equations is represented by the mathematical analysis of the non-linear
oscillations, which in many situations can be described by the equation \cite%
{Mic,Nay}
\begin{equation}
\ddot{x}+f\left(\dot{x}\right)g(x)+h(x)=0,
\end{equation}
which is a particular case of Eq.~(\ref{Lev}). where $f$, $g$ and $h$ are
arbitrary $C^1$ functions. From a physical point of view the Li\'{e}nard
equation can also be interpreted as the generalization of the equation of
motion of the damped oscillator,
\begin{equation}
\ddot{x}+\gamma \dot{x}+\omega ^2 x=0,
\end{equation}
with $\gamma =\mathrm{constant}$ and $\omega ^2=\mathrm{constant}$,
respectively \cite{Ben}, describing, from a physical point of view,
nonlinear or anharmonic oscillations \cite{Eu,Ha3}. A new criterion for
integrability of the Li\'{e}nard equation using an approach based on
nonlocal transformations, was obtained in \cite{Sil1}. Some previously known
criteria for integrability were also reobtained, and several new examples of
integrable Li\'{e}nard equations were given. A new family of the
Li\'enard-type equations which admits a non-standard autonomous Lagrangian
was found in \cite{Sil2}, and autonomous first integrals for each member
were obtained. Four new integrability criteria of a particular type of Li%
\'{e}nard type equations were obtained in \cite{Sil3} by studying the
connections between this family of Li\'enard--type equations and type III
Painlev\'e--Gambier equations. The results were illustrated by providing
examples of some integrable Li\'enard--type equations.

The connection between the linear harmonic oscillator
equation and some classes of second order nonlinear ordinary differential equations of
Li\'{e}nard and generalized Li\'{e}nard type, which physically describe important oscillator
systems, was investigated in \cite{Liang}. By using a method inspired by quantum mechanics, and which consist
on the deformation of the phase space coordinates of the harmonic oscillator, one can
generalize the equation of motion of the classical linear harmonic oscillator to several
classes of strongly non-linear differential equations. The first integrals, and a number
of exact solutions of the corresponding equations can then be explicitly obtained. The procedure can be further generalized to derive explicit general solutions of nonlinear second order differential equations unrelated to the harmonic oscillator.

An important relation can be established between the Li\'{e}nard type
equations, and the first kind first order Abel differential equation \cite%
{kamke, Pol},
\be
\frac{dy}{dx}=p(x)y^3+q(x)y^2.
\ee

This relation allows to find
some exact general solutions of the Li\'{e}nard type equations by using the
integrability conditions of the Abel equation. Some general integrability
conditions for the Abel equation were obtained in \cite{Mak1,Mak2,Mak3}.

By using the Li\'{e}nard equation-Abel equation equivalence, a class of
exact solutions of the Li\'{e}nard differential equation was obtained in
\cite{Ha4}. With the use of an exact integrability condition for the Abel
equation (the Chiellini lemma) \cite{Chiel, Ha5}, as a first step one can
obtain the exact general solution of the Abel equation. Once the solution of
the Abel equation is known, a class of exact solutions of the Li\'{e}nard
equation, expressed in a parametric form can be found. The Chiellini
integrability condition of the Abel equation was extended to the case of the
generalized Abel equation in \cite{Ha5}. The exact solutions of some
particular Li\'{e}nard type equations, including a generalized van der Pol
type equation, were also explicitly obtained. The interesting connection
between dissipative nonlinear second order differential equations and the
Abel equations with cubic and quadratic terms was considered in \cite{Rosu1}%
. The general solution of the corresponding integrable dissipative equations
was derived by using the Chiellini integrability condition. In \cite{Rosu2}
nonsingular parametric oscillators, Darboux related to the classical
harmonic oscillator, and having periodic dissipative/gain properties, have
been found through a modified factorization method. The same method can be
also applied to the so-called upside-down (hyperbolic) "oscillator". The barotropic Friedmann-Robertson-Walker cosmologies in the comoving time lead in the radiation-dominated case to scale factors of identical form as for the Chiellini dissipative scale factors in conformal time \cite{Rosu3}.  This is due to the Ermakov equation,  which is obtained in this case. Some
results in the integrability of the Abel equations were discussed and
reformulated in \cite{Man4}. Analytic techniques for the solutions of
nonlinear oscillators with damping using the Abel equation, the use of the
Chiellini integrability condition for the study of planar isochronous
systems and of the Hamiltonian structures of the Li\'{e}nard equation, as
well as the relation between the generalized damped Milne-Pinney equation
and the Chiellini method were investigated in \cite{new1,new2,new3}.
Travelling waves solutions in reaction-diffusion-convection systems were
obtained, with the use of the Chiellini and Lemke integrability conditions,
in \cite{Har,Har1}.

There are several methods that can be used for obtaining an
exact solution of the general Abel equation
\begin{equation}  \label{Agen}
\frac{dy}{dx}=p(x)y^3 + q(x)y^2 + r(x)y + s(x).
\end{equation}

An important property of the Abel differential equations is that they can be systematized into equivalence classes \cite{Liouville, Cheb, Bouqet}.
Two Abel ordinary differential equations are defined to be equivalent if they can be obtained one from the other by means of the transformation \cite{Bouqet}
\be
\left(y,x\right)\rightarrow \left(Y,X\right):x=F(X),y=P(X)Y(X)+Q(X), F'P\neq 0,
\ee
where $X$ and $Y(X)$ are the new independent and dependent variables, while $F$, $P$ and $Q$ are arbitrary
functions of $X$. Then the equivalence class for an ordinary differential equation is the set of all ordinary differential equations equivalent to the given one. As shown in \cite{Liouville, Cheb, Bouqet}, to each equivalence class we can associate an infinite sequence of absolute invariants (see \cite{Cheb} for the definition of absolute and relative invariants). As shown initially in  \cite{Liouville} (see also \cite{Cheb} and \cite{Bouqet} for detailed discussions), one can associate to the Abel equation the relative invariant $S_3$ of weight 3,
\be
S_3(x)\equiv s(x)p^3(x)+\frac{1}{3}\left[\frac{2q^2(x)}{9}-r(x)q(x)p(x)+p(x)q'(x)-q(x)p'(x)\right].
\ee
$S_3$ can be used to recursively generate an infinite sequence of relative invariants $S_{2m+1}$ of weights $2m+1$  through the relations
\be
S_{2m+1}(x)=p(x)S'_{2m-1}(x)-(2m-1)S_{2m-1}(x)\left[p'(x)+r(x)p(x)-\frac{q^2(x)}{3}\right].
\ee
With the help of the relative invariants $S_{2m+1}$ one can construct the absolute invariants $I_n$, given by \cite{Cheb, Bouqet}
\be
I_1(x)=\frac{S_5^3(x)}{S_3^5(x)}, I_2(x)=\frac{S_3(x)S_7(x)}{S_5^2(x)}, I_3(x)=\frac{S_9(x)}{S_3^2(x)},....
\ee
If $I_1(x)$ is a constant, then all the other invariants are also constant. The early results on the invariants of the Abel equation were used in \cite{Cheb} to obtain a set of new integrable Abel ordinary differential equation classes, with some depending on arbitrary parameters, and to introduce
an explicit method for  verifying or refuting the equivalence between two given Abel ordinary differential equations. The same approach was also used in \cite{Bouqet} for the analytic study of the buoyancy-drag equation with a time-dependent acceleration $\gamma (t)$, by determining its equivalence class under the point transformations, which allowed to define for some values of $\gamma (t)$ a time-dependent Hamiltonian from which the buoyancy-drag equation can be derived.

There are also several other methods for obtaining solutions of the Abel equation. If we introduce the transformations  \cite{kamke, Murphy,Pol}
\begin{equation}
y(x)=\omega (x)\eta \left(\xi (x)\right)-\frac{q(x)}{3p(x)},
\end{equation}
with
\begin{equation}
\omega (x)=e^{\int{\left[r(x)-\frac{q^2(x)}{3p(x)}\right]dx}}, \xi(x)=\int{%
p(x)\omega ^2(x)dx},
\end{equation}
then the Abel equation can be reduced to the normal form
\begin{equation}
\frac{d\eta}{dx}=\eta ^3+I(x),
\end{equation}
where
\begin{equation}
p(x)\omega ^3(x)I(x)=s(x)+\frac{1}{3}\frac{d}{dx}\frac{q(x)}{p(x)}-\frac{%
r(x)q(x)}{3p(x)}+\frac{2q^3(x)}{27p^2(x)}.
\end{equation}
Then, if $I(x)$ is constant, the Abel equation can be integrated \cite{Man4}.

On the other hand, if $y = y_1(x)$ is a particular solution Eq.~(\ref{Agen}%
), then by means of the transformations \cite{kamke}
\begin{equation}
u(x) = \frac{E(x)}{y(x)-y_1(x)},
\end{equation}
where
\begin{equation}
E(x) = \exp \Bigg\{\int{\left[3p(x)y_1 ^2 + 2q(x)y_1 + r(x)\right] dx} %
\Bigg\},
\end{equation}
then the Abel equation can be transformed into the form
\begin{equation}
\frac{du}{dx}+\frac{\Phi _1}{u}+\Phi _2=0,
\end{equation}
where
\begin{equation}
\Phi _1 x)=p(x)E^2(x),\Phi _2(x) = \left[3p(x)y_(x) + q(x)\right]E(x).
\end{equation}
Therefore, if the particular solution is given by
\begin{equation}
y_1=-\frac{q(x)}{ 3p(x)},
\end{equation}
then $\Phi_2 = 0$, and the general solution of the Abel equation can be
found by integrating an ordinary differential with separable variables \cite%
{kamke}.

It is the goal of the present paper to present some {\it exact integrability
cases} of the {\it extended Li\'{e}nard equation}
\begin{equation}
y^{\prime \prime }+f\left(
y\right) \left(y^{\prime }\right)^{n}+k\left( y\right) \left(y^{\prime
}\right)^{m}+g\left(y\right) y^{\prime }+h\left( y\right) =0,
\end{equation}
where $f(y)$, $k(y)$, $g(y)$, and $h(y)$ are arbitrary functions. This second order highly
nonlinear differential equation can be transformed to an equivalent \textit{%
generalized first kind first order Abel type equation}.

As a first step in our analysis we study in detail the quadratic-cubic
extended Li\'{e}nard equations, with $n=2$ and $m=3$, respectively. For this
case, we first obtain the general solutions of the Li\'{e}nard equation
under the assumption that one \textit{particular solution of the associated
Abel equation} is known. Under this assumption, we obtain three
integrability cases of the extended quadratic-cubic Li\'{e}nard equation,
corresponding to $n=2$ and $m=3$.

With the use of the \textit{Chiellini
integrability condition} \cite{Chiel}, the general solutions of the Abel and
quadratic-cubic Li\'{e}nard equations with coefficients satisfying some
differential conditions can be obtained in an exact closed form. We also
show that if a particular solution of the Abel equation is known, the
general solution of the extended quadratic cubic Li\'{e}nard equation can be
obtained by quadratures, if the coefficients of the equation and the
particular solution satisfy a given differential condition.

As a next step in our study. we extend \textit{the Chiellini integrability
condition to generalized Abel equations} with $g(y)\equiv 0$ and $h(y)\equiv
0$, and arbitrary $n$ and $m$. Thus, with the help of the generalized
Chiellini Lemma, we obtain the general solution of the corresponding
extended Li\'{e}nard equation for arbitrary $m$ and $n$. The application of
the generalized Chiellini condition to the case of the reduced Riccati
equation is also considered, thus leading to a new integrability case of
this equation (see \cite{Ha,Ha1,Ha2} for alternative integrability
conditions of the Riccati equation).

In the present paper we use the term \textit{integrable or exactly
integrable differential equation} as corresponding to a differential
equation that can be solved by \textit{quadratures, or by integrals,
definite or indefinite}, respectively.

The present paper is organized as follows. In Section~\ref{sect1-1} we
define the quadratic-cubic Li\'{e}nard equation, and we obtain its
associated Abel equation. Then, by assuming that a particular solution of
the associated Abel equation is known, we present two cases of exact
integrability of the Abel and Li\'{e}nard equations, respectively. In Section~\ref{sect2a} we obtain the integrability condition of the extended Li\'{e}nard  equation under the assumption that the coefficients $g(y)$ and $h(y)$  identically vanish, so that $g(y)=h(y)\equiv 0$. We
discuss and conclude our results in Section~\ref{sect2}.

\section{Exact solutions of the extended Li\'{e}nard equation}\label{sect1-1}

In the following we denote $y^{\prime }=\frac{dy}{dx}$ and $y^{\prime \prime
}=\frac{d^2y}{dx^2}$, respectively. We introduce the extended Li\'{e}nard
equation by means of the following

\textbf{Definition}. {\it The second order non-linear ordinary differential
equation
\begin{equation}
y^{\prime \prime }+f\left( y\right) \left( y^{\prime }\right) ^{n}+k\left(
y\right) \left( y^{\prime }\right) ^{m}+g\left( y\right) y^{\prime }+h\left(
y\right) =0,  \label{1_1}
\end{equation}%
where $f(y),k(y),g(y),h(y)\in C^{\infty }(I)$ are arbitrary functions
defined on a real interval $I\subseteq \Re $, $f(y),k(y),g(y),h(y)\neq
0,\forall y\in I$, and $n,m\in \Re $ are constants satisfying the conditions
$m,n>0$, and $n\neq m$, respectively, is called the extended Li\'{e}nard
differential equation.}

By introducing a new dependent variable $u=y^{\prime }$, Eq.~(\ref{1_1})
becomes
\begin{equation}
\frac{du}{dy}+f\left( y\right) u^{n-1}+k\left( y\right) u^{m-1}+g\left(
y\right) +\frac{h\left( y\right) }{u}=0.  \label{2_1}
\end{equation}

Let $u=\frac{1}{v}$, then we rewrite Eq. (\ref{2_1}) as%
\begin{equation}\label{3}
\frac{dv}{dy}=f\left( y\right) v^{\alpha}+k\left( y\right) v^{\beta}+g\left(
y\right) v^{2}+h\left( y\right) v^{3}, \alpha =3-n, \beta =3-m.
\end{equation}

{\bf Definition} {\it We call Eq.~(\ref{3}) the generalized Abel equation associated to
the extended Li\'{e}nard equation}.

Hence once the solutions of the first
order differential equation (\ref{3}) are known, the solutions of the
corresponding extended Li\'{e}nard equation can also be found. In the
following we will obtain three classes of exact solutions of Eq.~(\ref{3})
for $n=2$ and $m=3$, which will generate three classes of exact solutions of
the extended Li\'{e}nard equation.

\subsection{The quadratic-cubic extended Li\'{e}nard equation}

Assume that $n=2$ and $m=3$, then the extended Li\'{e}nard equation takes
the form
\begin{equation}  \label{ls}
y^{\prime \prime }+f\left( y\right) \left(y^{\prime }\right)^{2}+k\left(
y\right) \left(y^{\prime }\right)^{3}+g\left( y\right) y^{\prime }+h\left(
y\right) =0.
\end{equation}
We call Eq.~(\ref{ls}) \textit{the quadratic-cubic extended Li\'{e}nard
equation}. For this choice of $n$ and $m$ the generalized associated Abel
Eq.~(\ref{3}) becomes a standard first kind first order Abel type
differential equation, given by%
\begin{equation}
\frac{dv}{dy}=k\left( y\right) +f\left( y\right) v+g\left( y\right)
v^{2}+h\left( y\right) v^{3}.  \label{4_1}
\end{equation}

We assume that a particular solution $v_{p}$ that satisfies Eq. (\ref{4_1})
is known, that is, $v_p$ satisfies the equation
\begin{equation}
\frac{dv_{p}}{dy}=k\left( y\right) +f\left( y\right) v_{p}+g\left( y\right)
v_{p}^{2}+h\left( y\right) v_{p}^{3}.  \label{b1}
\end{equation}%
Subtracting Eq.~(\ref{4_1}) and Eq.~(\ref{b1}) and introducing the new
function $F$ defined as $F=v-v_{p}$ then we obtain the Abel equation for $%
F\left( y\right) $%
\begin{equation}
\frac{dF}{dy}=\left[ f\left( y\right) +2g\left( y\right) v_{p}+3h\left(
y\right) v_{p}^{2}\right] F+\left[ g\left( y\right) +3h\left( y\right) v_{p}%
\right] F^{2}+h\left( y\right) F^{3}.  \label{b2}
\end{equation}

By introducing a new function $w$, defined as
\begin{equation}
U(y)=E(y)w(y),
\end{equation}%
where
\begin{equation}
E(y)=e^{\int \left[ f\left( y\right) +2g\left( y\right) v_{p}+3h\left(
y\right) v_{p}^{2}\right] dy},
\end{equation}%
it follows that Eq.~(\ref{b2}) takes the form
\begin{equation}
\frac{dw}{dy}=\left[ g\left( y\right) +3h\left( y\right) v_{p}\right]
Ew^{2}+h\left( y\right) E^{2}w^{3},  \label{b6}
\end{equation}%
or, equivalently,
\begin{equation}
\frac{dw}{dy}=\left[ g\left( y\right) +3h\left( y\right) v_{p}\right]
e^{\int \left[ f\left( y\right) +2g\left( y\right) v_{p}+3h\left( y\right)
v_{p}^{2}\right] dy}w^{2}+h\left( y\right) e^{2\int \left[ f\left( y\right)
+2g\left( y\right) v_{p}+3h\left( y\right) v_{p}^{2}\right] dy}w^{3}.
\label{b9}
\end{equation}

\subsubsection{General solution of the Abel and extended Li\'{e}nard
equation for $v_p=v_{p0}=-g(y)/3h(y)$}

If the particular solution of the Abel equation (\ref{4_1}) is given by
\begin{equation}
v_{p0}(y)=-\frac{g(y)}{3h(y)},  \label{ps}
\end{equation}%
then Eq.~(\ref{b6}) becomes
\begin{equation}
\frac{dw(y)}{dy}=h\left( y\right) E^{2}(y)w^{3}(y),E(y)=e^{\int {\left[
f(y)-g^{2}(y)/3h(y)\right] dy}}.  \label{24}
\end{equation}%
Eq.~(\ref{24}) can be straightforwardly integrated, giving
\begin{equation}
w(y)=\pm \frac{1}{\sqrt{2}\sqrt{C-\int h(y)e^{2\int {\left[
f(y)-g^{2}(y)/3h(y)\right] dy}}dy}},
\end{equation}%
where $C$ is an arbitrary constant of integration.

Then we immediately obtain
\begin{equation}
U(y)=\pm \frac{e^{\int {\left[ f(y)-g^{2}(y)/3h(y)\right] dy}}}{\sqrt{2}%
\sqrt{C-\int h(y)e^{2\int {\left[ f(y)-g^{2}(y)/3h(y)\right] dy}}dy}}.
\end{equation}

Hence
\begin{equation}
v=\frac{1}{y^{\prime }}=U(y)+v_{p0}(y)=\pm \frac{e^{\int {\left[
f(y)-g^{2}(y)/3h(y)\right] dy}}}{\sqrt{2}\sqrt{C-\int h(y)e^{2\int {\left[
f(y)-g^{2}(y)/3h(y)\right] dy}}dy}}-\frac{g(y)}{3h(y)}.
\end{equation}

Therefore we have obtained the following

\textbf{Theorem 1}. {\it  If the generalized Abel equation associated to the
extended quadratic-cubic Li\'{e}nard equation with $n=2$ and $m=3$ has the
particular solution given by Eq.~(\ref{ps}), and the coefficients of the Li%
\'{e}nard equation satisfy the condition
\begin{equation}
\frac{d}{dy}\left[ \frac{g(y)}{h(y)}\right] =-3k(y)+\frac{f(y)g(y)}{h(y)}-%
\frac{2}{9}\frac{g^{3}(y)}{h^{2}(y)},  \label{31}
\end{equation}%
then the general solution of the quadratic-cubic Li\'{e}nard equation is
given by
\begin{equation}
x-x_{0}=\int \left\{ \pm \frac{e^{\int {\left[ f(y)-g^{2}(y)/3h(y)\right] dy}%
}}{\sqrt{2}\sqrt{C-\int h(y)e^{2\int {\left[ f(y)-g^{2}(y)/3h(y)\right] dy}%
}dy}}-\frac{g(y)}{3h(y)}\right\} dy,
\end{equation}%
where $x_{0}$ is an arbitrary integration constant. Eq.~(\ref{31}) is an
Abel equation for $g(y)$, and a Riccati equation for $1/h(y)$.}

\subsubsection{General solutions of Abel and extended Li\'{e}nard equation
for $v_p \neq -g(y)/3h(y)$ via the Chiellini integrability condition}

If the particular solution of the generalized Abel equation associated to
the extended Li\'{e}nard equation is not given by Eq.~(\ref{ps}), $v_p \neq
-g(y)/3h(y)$, a solution of Eq.~(\ref{b9}) can be found with the use of the
Chiellini Lemma, which can be formulated as follows \cite{Chiel,kamke}

\textbf{Lemma 1 (Chiellini (1931))} \cite{Chiel} {\it A first kind Abel type
differential equation of the form $y^{\prime }+p(x)y^3+q(x)y^{2}=0$ can be exactly
integrated if the functions $q(x)$ and $p(x)$ satisfy the condition
\begin{equation}
\frac{d}{dx}\left[ \frac{p(x)}{q(x)}\right] =Sq(x),\qquad S=\mathrm{constant}%
,\qquad S\neq 0.
\end{equation}
}

In order to solve the Abel equation with the help of the \textbf{Chiellini
Lemma} we introduce a new function $\theta (x)$, defined as
\begin{equation}
y(x)=\frac{q(x)}{p(x)}\theta (x).
\end{equation}

Then, with the use of the Chiellini integrability condition the Abel
equation can be transformed into a first order separable differential
equation,
\begin{equation}
\frac{d\theta }{dx}=\frac{q^{2}(x)}{p(x)}\theta \left( \theta ^{2}+\theta
+S\right) ,
\end{equation}%
with the general solution given by
\begin{equation}
\frac{p(x)}{q(x)}=K_{0}^{-1}e^{G(\theta ,S)},  \label{sol}
\end{equation}%
where $K_{0}^{-1}$ is an arbitrary constant of integration, and
\begin{equation}
e^{G(\theta ,S)}=\left\{
\begin{array}{lll}
\frac{\theta }{\sqrt{\theta ^{2}+\theta +S}}\exp \left[ -\frac{1}{\sqrt{4S-1}%
}\arctan \left( \frac{1+2\theta }{\sqrt{4S-1}}\right) \right] ,\qquad S>%
\frac{1}{4}, &  &  \\
&  &  \\
\frac{\theta }{1+2\theta }e^{1/(1+2\theta )},\qquad S=\frac{1}{4}, &  &  \\
&  &  \\
\frac{\theta }{\sqrt{\theta ^{2}+\theta +S}}\exp \left[ \frac{1}{\sqrt{1-4S}}%
\mathrm{arctanh}\left( \frac{1+2\theta }{\sqrt{1-4S}}\right) \right] ,\qquad
S<\frac{1}{4}, &  &
\end{array}%
\right.  \label{sol2}
\end{equation}%
respectively. Eq.~(\ref{sol}) determines $\theta $ as a function of $x$ and
of the constant $S$, respectively.

For the Abel Eq.~(\ref{b9}) the Chiellini integrability condition is given
by
\begin{equation}
\frac{d}{dy}\left\{ \frac{h\left( y\right) e^{\int \left[ f\left( y\right)
+2g\left( y\right) v_{p}+3h\left( y\right) v_{p}^{2}\right] dy}}{\left[
g\left( y\right) +3h\left( y\right) v_{p}\right] }\right\} =S\left[ g\left(
y\right) +3h\left( y\right) v_{p}\right] e^{\int \left[ f\left( y\right)
+2g\left( y\right) v_{p}+3h\left( y\right) v_{p}^{2}\right] dy},  \label{b10}
\end{equation}%
or, equivalently
\begin{equation}
\frac{d}{dy}\left[ \frac{E(y)}{v_{p}(y)-v_{p0}(y)}\right] =9Sh(y)\left[
v_{p}(y)-v_{p0}(y)\right] E(y).  \label{intc}
\end{equation}

\subsubsection{The case $E(y)\equiv 1$}

If the particular solution of the Abel equation (\ref{4_1}) satisfies the
condition
\begin{equation}
3h(y)v_p^2+2g(y)v_p+f(y)=0,
\end{equation}
or
\begin{equation}  \label{psg}
v_p(y)=-\frac{g(y)}{3h(y)}\pm\sqrt{\frac{g^2(y)}{9h^2(y)}-\frac{f(y)}{3h(y)}}%
,
\end{equation}
then, without any loss of generality, one can choose $E(y)\equiv 1$. For
this case the Chiellini integrability condition becomes
\begin{equation}
\frac{d}{dy}\left[ \frac{1}{v_{p}(y)-v_{p0}(y)}\right] =9Sh(y)\left[
v_{p}(y)-v_{p0}(y)\right] ,
\end{equation}%
or
\begin{equation}
\frac{d}{dy}\left\{ \left[ \frac{g^{2}(y)}{9h^{2}(y)}-\frac{f(y)}{3h(y)}%
\right] ^{-1/2}\right\} =9Sh(y)\sqrt{\frac{g^{2}(y)}{9h^{2}(y)}-\frac{f(y)}{%
3h(y)}}.  \label{psg2}
\end{equation}

By taking into account the identity
\be
\frac{1}{\sqrt{b(y)}}\frac{d}{dy}\frac{1%
}{\sqrt{b(y)}}=-\frac{1}{2b^{2}(y)}\frac{d}{dy}b(y),
\ee
the condition (\ref{psg2}) can be first reformulated as
\begin{equation}
\frac{d}{dy}\left[ \frac{g^{2}(y)}{9h^{2}(y)}-\frac{f(y)}{3h(y)}\right]
=-18Sh(y)\left[ \frac{g^{2}(y)}{9h^{2}(y)}-\frac{f(y)}{3h(y)}\right] ^{2}.
\end{equation}

Hence we have obtained the following

\textbf{Theorem 2}.{\it  a) If the coefficients of the general Abel equation (\ref%
{4_1}) satisfy the condition
\begin{equation}
\frac{g^{2}(y)}{3h^{2}(y)}-\frac{f(y)}{h(y)}=\frac{1}{6S\int h(y)dy+C_{0}/3},
\label{37}
\end{equation}%
where $C_{0}$ is an arbitrary constant, then the Abel equation is exactly
integrable.

b) The general solution of the Abel equation (\ref{4_1}) with coefficients
satisfying condition (\ref{37}) is given by
\begin{equation}
v(y,S)=\pm \frac{\sqrt{g^{2}(y)-3f(y)h(y)}}{h(y)}\theta (y,S),
\end{equation}%
where $\theta (y,S)$ is a solution of the algebraic equation
\begin{equation}
K_{0}^{-1}e^{G\left[ \theta (y),S\right] }=\pm \frac{h(y)}{\sqrt{%
g^{2}(y)-3f(y)h(y)}}.
\end{equation}%
c) The extended quadratic-cubic Li\'{e}nard equation (\ref{ls}) with
coefficients satisfying condition (\ref{37}) is exactly integrable, and its
general solution is given by
\begin{equation}
x-x_{0}=\pm \frac{1}{3}\int \frac{1}{h(y)}\left\{ \sqrt{g^{2}(y)-3f(y)h(y)}%
\left[ 3\theta \left( y,S\right) +1\right] -g(y)\right\} dy,
\end{equation}%
where $x_{0}$ is an arbitrary constant of integration.}

\subsection{The general solution of the extended quadratic-cubic Li\'{e}%
nard for arbitrary particular solutions of the associated Abel equation}

We consider now the case of an arbitrary particular solution of the
associated Abel equation. From the Chiellini integrability condition (\ref%
{b10}), it follows that the arbitrary coefficient $f\left( y\right) $ of the
quadratic-cubic Li\'{e}nard equation must satisfy the following differential
condition
\begin{equation}
f\left( y\right) =\frac{d^{2}}{dy^{2}}\left[ \ln \left\vert \frac{g\left(
y\right) }{h\left( y\right) }+3v_{p}\right\vert \right] +S\frac{d}{dy}%
\left\{ h\left( y\right) \left[ \frac{g\left( y\right) }{h\left( y\right) }%
+3v_{p}\right] ^{2}\right\} -2g\left( y\right) v_{p}-3h\left( y\right)
v_{p}^{2}.  \label{b11}
\end{equation}%
With the help of the transformation given by%
\begin{equation}
w=\frac{\left[ g\left( y\right) +3h\left( y\right) v_{p}\right] }{h\left(
y\right) e^{\int \left[ f\left( y\right) +2g\left( y\right) v_{p}+3h\left(
y\right) v_{p}^{2}\right] dy}}\theta ,  \label{b12}
\end{equation}%
immediately inserting Eq. (\ref{b12}) into Eq. (\ref{b9}), the latter
becomes a separate variable differential equation

\begin{equation}
\frac{d\theta }{dy}=\frac{\left[ g\left( y\right) +3h\left( y\right) v_{p}%
\right] ^{2}}{h\left( y\right) }\theta \left( \theta ^{2}+\theta +S\right) .
\label{b13}
\end{equation}%
Or equivalently,
\begin{equation}
\int \frac{1}{\theta \left( \theta ^{2}+\theta +S\right) }d\theta =\int
\frac{\left[ g\left( y\right) +3h\left( y\right) v_{p}\right] ^{2}}{h\left(
y\right) }dy.  \label{b14}
\end{equation}%
Now in view of the relation
\begin{equation}
w=\frac{\left[ g\left( y\right) +3h\left( y\right) v_{p}\right] }{h\left(
y\right) e^{\int \left[ f\left( y\right) +2g\left( y\right) v_{p}+3h\left(
y\right) v_{p}^{2}\right] dy}}\theta =\left( v-v_{p}\right) e^{-\int \left[
f\left( y\right) +2g\left( y\right) v_{p}+3h\left( y\right) v_{p}^{2}\right]
dy},
\end{equation}%
we obtain%
\begin{equation}
v\left( y\right) =\frac{g\left( y\right) +3h\left( y\right) v_{p}}{h\left(
y\right) }\theta \left( y\right) +v_{p}.  \label{b15}
\end{equation}%
Next substituting Eq. (\ref{b11}) into Eq. (\ref{4_1}), the latter takes the
form
\begin{equation}
\frac{dv}{dy}=k\left( y\right) +\Bigg\{\frac{d^{2}}{dy^{2}}\left[ \ln
\left\vert \frac{g\left( y\right) }{h\left( y\right) }+3v_{p}\right\vert %
\right] + \\
S\frac{d}{dy}\left\{ h\left( y\right) \left[ \frac{g\left( y\right) }{%
h\left( y\right) }+3v_{p}\right] ^{2}\right\} -2g\left( y\right)
v_{p}-3h\left( y\right) v_{p}^{2}\Bigg\}v+g\left( y\right) v^{2}+h\left(
y\right) v^{3}.  \label{b16}
\end{equation}

Therefore we have obtained the following

\textbf{Theorem 3.} {\it a) If a particular solution $v_p$ of the Abel Eq.~(\ref%
{4_1}) is known, then the general solution of the Abel Eq.~(\ref{b16}) is
given by Eq.~(\ref{b15}).

b) The general solution of the quadratic-cubic extended Li\'{e}nard equation
with coefficient $f(y)$ satisfying the condition (\ref{b11}) is given by
\begin{equation}
x-x_{0}=\int \left[ \frac{g\left( y\right) +3h\left( y\right) v_{p}}{h\left(
y\right) }\theta \left( y\right) +v_{p}\right] dy,
\end{equation}%
where $x_{0}$ is the arbitrary integration constant.}

The function $\theta \left( y\right) $, given by Eq. (\ref{b15}), and which
determines the solution of the Abel equation, can be determined from Eq. (%
\ref{b14}). Subsequently, we have obtained the general solution of the
extended quadratic-cubic Li\'{e}nard Eq.~(\ref{ls}) under the assumption
that a particular solution of the associated Abel equation is known.

\section{General solution of the extended Li\'{e}nard equation for
arbitrary $n$, $m$, and $g(y)=h(y)\equiv 0$}\label{sect2a}

We assume now that the coefficients $g(y)$ and $h(y)$ of the extended Li\'{e}%
nard equation identically vanish, so that $g(y)=h(y)\equiv 0$. In this case
the extended Li\'{e}nard equation and its associated Abel equation take the
form
\begin{equation}  \label{lg}
y^{\prime \prime }+f\left( y\right) \left( y^{\prime }\right) ^{n}+k\left(
y\right) \left( y^{\prime }\right) ^{m} =0,
\end{equation}
and
\begin{equation}  \label{ag}
\frac{dv}{dy}=f\left( y\right) v^{\alpha}+k\left( y\right) v^{\beta}, \alpha
=3-n, \beta =3-m,
\end{equation}
respectively, where $v=1/y^{\prime }$. Then we can formulate the following
\textit{generalized Chiellini} Lemma as

\textbf{Lemma 2 (Generalized Chiellini Lemma)}. {\it  If the coefficients $f(y)$
and $k(y)$ of the generalized Abel equation satisfy the differential
conditions
\begin{equation}
\frac{d}{dy}\left[ \left( \frac{f(y)}{k(y)}\right) ^{\frac{1}{\beta -\alpha}}%
\right]=SP^{\beta -1}f(y)\left[ \frac{f(y)}{k(y)}\right] ^{\frac{\alpha}{%
\beta -\alpha}}, \alpha \neq \beta,
\end{equation}
or
\begin{equation}
\frac{d}{dy}\left[ \left( \frac{f(y)}{k(y)}\right) ^{\frac{1}{\beta -\alpha}}%
\right]=SP^{\beta -1}k(y)\left[ \frac{f(y)}{k(y)}\right] ^{\frac{\beta }{%
\beta - \alpha}} , \alpha \neq \beta,
\end{equation}
where $\alpha$, $\beta $, $S$ and $P$ are arbitrary constants, then the
generalized Abel equation is exactly integrable.}

\textbf{Proof.} By means of a simple transformation the generalized Abel
Eq.~(\ref{ag}) can be transformed into the equivalent form
\begin{equation}
\frac{d\ln v}{dy}=f\left( y\right) v^{\alpha -1}+k\left( y\right) v^{\beta
-1}.
\end{equation}

Let's assume now the function $v$ can be represented as
\begin{equation}
v(y)=F\left[ f(y),k(y)\right] \theta (y),  \label{trans}
\end{equation}%
where $F$ and $\theta $ are two arbitrary functions. Then we have $\ln v=\ln
F+\ln \theta $, and
\begin{equation}
\frac{d\ln v}{dy}=\frac{1}{F}\frac{dF}{dy}+\frac{1}{\theta }\frac{d\theta }{%
dy}=fF^{\alpha -1}\theta ^{\alpha -1}+kF^{\beta -1}\theta ^{\beta -1}.
\label{63}
\end{equation}%
Let's assume that the condition
\begin{equation}
\frac{1}{S}\frac{1}{F}\frac{dF}{dy}=P^{\beta -\alpha }fF^{\alpha
-1}=kF^{\beta -1},
\end{equation}%
where $S$ and $P$ are arbitrary constants, holds for all $F$, $f$ and $k$.
Then the above condition determines first the function $F$ as
\begin{equation}
F=P\left( \frac{f}{k}\right) ^{\frac{1}{\beta -\alpha }}.
\end{equation}%
Then it follows that the function $F$ must satisfy the conditions
\begin{equation}
\frac{d}{dy}\left[ \left( \frac{f}{k}\right) ^{\frac{1}{\beta -\alpha }}%
\right] =SP^{\beta -1}f\left( \frac{f}{k}\right) ^{\frac{\alpha }{\beta
-\alpha }},
\end{equation}%
or, equivalently,
\begin{equation}
\frac{d}{dy}\left[ \left( \frac{f}{k}\right) ^{\frac{1}{\beta -\alpha }}%
\right] =SP^{\beta -1}k\left( \frac{f}{k}\right) ^{\frac{\beta }{\beta
-\alpha }}.
\end{equation}%
Hence Eq.~(\ref{trans}) becomes
\begin{equation}
v=F\left( f,k\right) \theta =P\left( \frac{f}{k}\right) ^{\frac{1}{\beta
-\alpha }}\theta .
\end{equation}

Eq.~(\ref{63}) immediately gives
\begin{equation}
SkF^{\beta -1}+\frac{1}{\theta }\frac{d\theta }{dy}=P^{\alpha -\beta
}kF^{\beta -1}\theta ^{\alpha-1}+kF^{\beta -1}\theta ^{\beta -1},
\end{equation}
which leads to the following equation with separable variables for $\theta $%
,
\begin{equation}
\frac{d\theta (y) }{dy}=P^{\beta -1}k(y)\left[\frac{f(y)}{k(y)}\right] ^{%
\frac{\beta -1}{\beta -\alpha}}\left[P^{\alpha -\beta }\theta ^{\alpha
}(y)+\theta ^{\beta }(y)-S\theta (y)\right].
\end{equation}
This ends the proof of the \textit{Generalized Chiellini Lemma}.

By using \textbf{Lemma 2} we can obtain the general solution of the extended
Li\'{e}nard equation (\ref{lg}) by means of the following

\textbf{Theorem 4}. {\it If the coefficients $f(y)$ and $k(y)$ of the extended Li%
\'{e}nard equation (\ref{lg}) satisfy the condition
\begin{equation}
\frac{d}{dy}\left[ \left( \frac{f(y)}{k(y)}\right) ^{\frac{1}{n-m}}\right]
=SP^{2-m}f(y)\left[ \frac{f(y)}{k(y)}\right] ^{\frac{3-n}{n-m}},n\neq m,
\end{equation}%
or
\begin{equation}
\frac{d}{dy}\left[ \left( \frac{f(y)}{k(y)}\right) ^{\frac{1}{n-m}}\right]
=SP^{2-m}k(y)\left[ \frac{f(y)}{k(y)}\right] ^{\frac{3-m}{n-m}},n\neq m,
\end{equation}%
where $S$ and $P$ are arbitrary constants, then the general solution of the
extended Li\'{e}nard equation can be obtained as
\begin{equation}
x-x_{0}=P\int {\left[ \frac{f(y)}{k(y)}\right] ^{\frac{1}{n-m}}\theta (y)dy},
\end{equation}%
where $x_{0}$ is the arbitrary integration constant and $\theta (y)$ is the
solution of the equation
\begin{equation}
\theta (y)=H^{-1}\left\{ P^{2-m}\int {k(y)\left[ \frac{f(y)}{k(y)}\right] ^{%
\frac{2-m}{n-m}}dy}\right\} ,n\neq m.  \label{bbb}
\end{equation}%
and we have denoted
\begin{equation}
H(\theta )=\int {\frac{d\theta }{P^{m-n}\theta ^{3-n}+\theta ^{3-m}-S\theta }%
=P^{2-m}}\int {k(y)\left[ \frac{f(y)}{k(y)}\right] ^{\frac{2-m}{n-m}}dy}%
,n\neq m.  \label{aaa}
\end{equation}
}

The proof of \textbf{Theorem 4} is immediate.

\subsection{The Chiellini integrability condition for the Riccati equation}

In the case $\alpha =0$, $\beta =2$, the generalized Abel equation (\ref{ag}%
) becomes a reduced Riccati type equation, given by
\begin{equation}  \label{ric}
\frac{dv}{dy}=f\left( y\right)+k\left( y\right) v^2.
\end{equation}
Then the generalized Chiellini Lemma allows us to obtain an integrability
condition for the Riccati equation, which is given by the following

\textbf{Theorem 5.} {\it If the coefficients $f(y)$ and $g(y)$ of the reduced
Riccati equation (\ref{ric}) satisfy the condition
\begin{equation}
\frac{d}{dy}\left[ \sqrt{\frac{f(y)}{k(y)}}\right] =Kf(y),  \label{condr}
\end{equation}%
where $K$ is an arbitrary constant, then the Riccati equation is exactly
integrable, and its solution is given by
\begin{equation}
v(y)=\sqrt{\frac{f(y)}{k(y)}}\left\{ \sqrt{1-\frac{K^{2}}{4}}\tan \left\{
\sqrt{1-\frac{K^{2}}{4}}\left[ \int {\sqrt{f(y)k(y)}dy}+C\right] \right\} +%
\frac{K}{2}\right\} ,K\neq 2,
\end{equation}%
\begin{equation}
v(y)=\sqrt{\frac{f(y)}{k(y)}}\left[ -\frac{1}{\int {\sqrt{f(y)k(y)}dy}+C}\mp
1\right] ,K=\mp 2.
\end{equation}
}

\textbf{Proof.} We look for a solution of the Riccati equation (\ref{ric})
of the form
\begin{equation}\label{82}
v(y)=\sqrt{\frac{f(y)}{k(y)}}\theta (y).
\end{equation}%

By substituting into the Riccati equation we obtain
\begin{equation}
\frac{d}{dy}\left[ \sqrt{\frac{f(y)}{k(y)}}\right] \theta (y)+\sqrt{\frac{%
f(y)}{k(y)}}\frac{d\theta (y)}{dy}=f(y)\left[ 1+\theta ^{2}(y)\right] .
\label{79}
\end{equation}%

With the use of the condition~(\ref{condr}), Eq.~(\ref{79}) becomes a first
order separable differential equation, which can be integrated as
\begin{equation}
\int {\frac{d\theta }{\theta ^{2}-K\theta +1}}=\int {\sqrt{f(y)k(y)}dy}.
\end{equation}%
Hence we obtain
\begin{equation}
\theta (y)=\sqrt{1-\frac{K^{2}}{4}}\tan \left\{ \sqrt{1-\frac{K^{2}}{4}}%
\left[ \int {\sqrt{f(y)k(y)}dy}+C\right] \right\} +\frac{K}{2},K\neq 2,
\end{equation}%
\begin{equation}
\theta (y)=-\frac{1}{\int {\sqrt{f(y)k(y)}dy}+C}\mp 1,K=\mp 2,
\end{equation}
where $C$ is an arbitrary integration constant. Therefore the general
solution of the Riccati Eq.~(\ref{ric}) with coefficients satisfying a
generalized Chiellini type integrability immediately follows from Eq.~(\ref{82}), and this completes the Proof of the Theorem.

\section{Concluding remarks}\label{sect2}

In the present paper, by using the Chiellini integrability condition of the
Abel equation we have obtained three exact solutions of the extended
quadratic-cubic Li\'{e}nard equation, and of the associated Abel equation.
The solutions have been obtained under the assumption that a particular
solution of the Abel equation is known. We have considered the cases in
which these particular solutions have a specific form, as well as the
general case in which the functional form of $v_{p}$ was not imposed in
advance. Once the particular solutions of the associated Abel equation are
known, the general solution of the nonlinear quadratic-cubic second order
extended Li\'{e}nard differential equation can be obtained by quadratures,
if the four coefficients of the equation satisfy some consistency conditions.

We have also presented a generalization of the Chiellini integrability
condition for the case of generalized Abel equations of the type (\ref{ag}),
which gives the possibility of obtaining some exact solutions of the
extended Li\'{e}nard equation for arbitrary $n$ and $m$, and for $%
g(y)=h(y)\equiv 0$. As a particular application of the generalized Chiellini
integrability condition we have also obtained an integrability case of the
reduced Riccati equation.

All the solutions obtained in the present paper require the existence of
some differential relations between the coefficients of the differential
equations. These constraints impose strong restrictions on the functional
form of the solutions, thus limiting the number of possible solutions of the
Abel and extended Li\'{e}nard equations that can be obtained in this way.
Some applications of the present results to some extended Li\'{e}nard
equations that describe some physically interesting mathematical models will
be presented in a future paper.

\section*{Acknowledgments}

We would like to thank the anonymous referee for comments and suggestions that helped us to improve our manuscript.


\begin{thebibliography}{99}
\bibitem{Lien} A. Li\'{e}nard, \'{E}tude des oscillations entretenues, Revue
g\'{e}n\'{e}rale de l'\'{e}lectricit\'{e} \textbf{23}, 901-912 (1928).

\bibitem{Lien1} A. Li\'{e}nard, \'{E}tude des oscillations entretenues,
Revue g\'{e}n\'{e}rale de l'\'{e}lectricit\'{e} \textbf{23}, 946-954 (1928).

\bibitem{Lev} N. Levinson and O. Smith, A general equation for relaxation
oscillations, Duke Mathematical Journal \textbf{9}, 382-403 (1942).

\bibitem{2} A. A. Andronov, E. A. Leontovich, I. I. Gordon, and A. G. Maier,
Qualitative Theory of Second Order Dynamic Systems, New York, Wiley (1973).

\bibitem{lit0} V. K. Chandrasekar, M. Senthilvelan, A. Kundu, and M.
Lakshmanan, A nonlocal connection between certain linear and nonlinear
ordinary differential equations/oscillators, Journal of Physics A:
Mathematical and Theoretical \textbf{39}, 9743-9754 (2006).

\bibitem{lit1} X.-G. Liu, M.-L. Tang, and R. R. Martin, Periodic solutions
for a kind of Li\'{e}nard equation, Journal of Computational and Applied
Mathematics \textbf{219}, 263-275 (2008).

\bibitem{lit2} L. Zou, X. Chen, and W. Zhang, Local bifurcations of critical
periods for cubic Li\'{e}nard equations with cubic damping, Journal of
Computational and Applied Mathematics \textbf{222}, 404-410 (2008).

\bibitem{lit3} S. N. Pandey, P. S. Bindu, M. Senthilvelan, and M.
Lakshmanan, A group theoretical identification of integrable cases of the Li%
\'{e}nard-type equation $x"+f(x)x^{\prime }+g(x)=0$. I. Equations having
nonmaximal number of Lie point symmetries, Journal of Mathematical Physics
\textbf{50} 082702-082702-19 (2009).

\bibitem{lit4} S. N. Pandey, P. S. Bindu, M. Senthilvelan, and M.
Lakshmanan, A group theoretical identification of integrable equations in
the Li\'{e}nard-type equation $x"+f(x)x^{\prime}+g(x)=0$. II. Equations having maximal
Lie point symmetries, Journal of Mathematical Physics \textbf{50},
102701-102701-25 (2009).

\bibitem{lit5} D. Banerjee and J. K. Bhattacharjee, Renormalization group
and Li\'{e}nard systems of differential equations, Journal of Physics A:
Mathematical and Theoretical \textbf{43}, 062001 (2010).

\bibitem{lit6} M. Messias, G. Alves, and R. M\'{a}rcio, Time-periodic
perturbation of a Li\'{e}nard equation with an unbounded homoclinic loop,
Physica D: Nonlinear Phenomena \textbf{240}, 1402-1409 (2011).

\bibitem{Mic} R. E. Mickens, Truly Nonlinear Oscillations: Harmonic Balance,
Parameter Expansions, Iteration, and Averaging Methods, World Scientific,
New Jersey-London-Singapore-Beijing-Shanghai-Hong Kong-Taipei-Chennai (2010).

\bibitem{Nay} A. H. Nayfeh and D. T. Mook, Nonlinear Oscillations, John
Wiley \& Sons, New York, Chichester (1995).

\bibitem{Ben} E. DiBenedetto, Classical mechanics: theory and mathematical
modeling, New York, N. Y., Birkh$\ddot{\mathrm{a}}$user, Springer, (2011).

\bibitem{Eu} N. Euler, Transformation properties of $\frac{d^2x}{dt^2} +f_1
(t)\frac{dx}{dt} +f_2 (t) x+f_3 (t) x^n = 0$, Journal of Nonlinear
Mathematical Physics \textbf{4}, 310-337 (1997).

\bibitem{Ha3} T. Harko, F. S. N. Lobo, and M. K. Mak, Integrability cases
for the anharmonic oscillator equation, Journal of Pure and Applied
Mathematics: Advances and Applications \textbf{10}, 115-129 (2013).

\bibitem{Sil1} N. A. Kudryashov and D. I. Sinelshchikov, On the criteria for
integrability of the Li\'{e}nard equation, Applied Mathematics Letters
\textbf{57} 114 - 120 (2016).

\bibitem{Sil2} N. A. Kudryashov and D. I. Sinelshchikov, New non-standard
Lagrangians for the Li\'{e}nard-type equations, Applied Mathematics Letters
\textbf{63} 124 - 129 (2017).

\bibitem{Sil3} N. A. Kudryashov and D. I. Sinelshchikov, On the
integrability conditions for a family of the Li\'{e}nard-type equations,  Regular and Chaotic Dynamics  {\bf 21}, 548 - 555 (2016).

\bibitem{Liang} T. Harko and S.-D. Liang, Exact solutions of the Li\'{e}nard and generalized Li\'{e}nard type ordinary non-linear differential equations obtained by deforming the phase space coordinates of the linear harmonic oscillator, Journal of Engineering Mathematics {\bf 98}, 93-111 (2016).

\bibitem{kamke} Kamke E. Differentialgleichungen: L\"{o}sungsmethoden und L%
\"{o}sungen. Chelsea, New York (1959).

\bibitem{Murphy} G. M. Murphy, Ordinary Differential Equations and their Solutions, Van Nostrand, Princeton (1960).

\bibitem{Pol} Polyanin A. D. and Zaitsev V. F. Handbook of Exact Solutions
for Ordinary Differential Equations. Chapman \& Hall/CRC, Boca Raton,
London, New York, Washington, D. C. (2003).

\bibitem{Mak1} M. K. Mak, H. W. Chan, and T. Harko, Solutions generating
technique for Abel-type nonlinear ordinary differential equations, Computers
\& Mathematics with Applications \textbf{41}, 1395-1401 (2001).

\bibitem{Mak2} M. K. Mak and T. Harko, New method for generating general
solution of Abel differential equation, Computers \& Mathematics with
Applications \textbf{43}, 91-94 (2002).

\bibitem{Mak3} T. Harko and M. K. Mak, Relativistic dissipative cosmological
models and Abel differential equation, Computers \& Mathematics with
Applications \textbf{46}, 849-853 (2003).

\bibitem{Ha4} T. Harko, F. S. N. Lobo, and M. K. Mak, A class of exact
solutions of the Li\'{e}nard type ordinary non-linear differential equation,
Journal of Engineering Mathematics \textbf{89}, 193-205 (2014).

\bibitem{Chiel} A. Chiellini, Sull'integrazione dell'equazione differenziale
$y^{\prime} +Py^2+Qy^ 3 =0$, Bollettino della Unione Matematica Italiana
\textbf{10}, 301-307 (1931).

\bibitem{Ha5} T. Harko, F. S. N. Lobo, and M. K. Mak, A Chiellini type
integrability condition for the generalized first kind Abel differential
equation, Universal Journal of Applied Mathematics \textbf{1}, 101-104
(2013).

\bibitem{Rosu1} S. C. Mancas and H. C. Rosu, Integrable dissipative
nonlinear second order differential equations via factorizations and Abel
equations, Phys. Lett. \textbf{A 377}, 1234-1238 (2013).

\bibitem{Rosu2} H. C. Rosu, O. Cornejo-Perez, and P. Chen, Nonsingular
parametric oscillators Darboux-related to the classical harmonic oscillator,
Europhys. Lett. \textbf{100}, 60006 (2012).

\bibitem{Rosu3} H. C. Rosu, S. C. Mancas, and P. Chen, Barotropic FRW cosmologies with Chiellini damping in comoving time, Mod. Phys. Lett. {\bf A 30}, 1550100 (2015).

\bibitem{Man4} S. C. Mancas and H. C. Rosu, Integrable Abel equations and
Vein's Abel equation, Math. Meth. Appl. Sci. \textbf{39}, 1376-1387 (2016).

\bibitem{new1} A. Ghose-Choudhury and P. Guha, An analytic technique for the
solutions of nonlinear oscillators with damping using the Abel Equation, to
be published in Discontinuity, Nonlinearity, and Complexity,
arXiv:1608.02324 [nlin.SI] (2016).

\bibitem{new2} A. Ghose Choudhury and P. Guha, Chiellini integrability
condition, planar isochronous systems and Hamiltonian structures of Li\'{e}%
nard equation, Discrete \& Continuous Dynamical Systems,  {\bf 22},  2465-2478 (2017).

\bibitem{new3} S. Mukherjee, A. Ghose Choudhury, and P. Guha, Generalized
damped Milne-Pinney equation and Chiellini method, arXiv:1603.08747
[nlin.SI] (2016).

\bibitem{Har} T. Harko and M. K. Mak, Travelling wave solutions of the
reaction-diffusion mathematical model of glioblastoma growth: An Abel
equation based approach, Mathematical Biosciences and Engineering \textbf{12}%
, 41-69 (2015).

\bibitem{Har1} T. Harko and M. K. Mak, Exact travelling wave solutions of
non-linear reaction-convection-diffusion equations -- an Abel equation based
approach, J. Math. Phys. 56, 111501 (2015).

\bibitem{Liouville} R. Liouville, Sur une \'{e}quation diff\'{e}rentielle du premier ordre, Acta mathematica {\bf 26},  55–78 (1902).

\bibitem{Cheb} E. S. Cheb-Terrab and A. D. Roche, Abel ODEs: equivalence and integrable classes, Computer Physics
Communications {\bf 130}, 204–231 (2000).

\bibitem{Bouqet} S. E. Bouqet, R. Conte, V. Kelsch, and F. Louvet, Solutions of the buoyancy-drag equation with a time-dependent acceleration, Journal of Nonlinear Mathematical Physics {\bf  24}, Supplement 1,   3–17 (2017).

\bibitem{Scal} P. Scalizzi, Soluzione di alcune equazioni del tipo di Abel,
Atti Accad. Naz. Lincei, Ser. \textbf{5}, 60-64 (1917).

\bibitem{Ha} M. K. Mak and T. Harko, New integrability case for the Riccati
equation, Applied Mathematics and Computation \textbf{218}, 10974-10981
(2012).

\bibitem{Ha1} M. K. Mak and T. Harko, New further integrability cases for
the Riccati equation, Applied Mathematics and Computation \textbf{219},
7465-7471 (2013).

\bibitem{Ha2} T. Harko, F. S. N. Lobo, and M. K. Mak, Analytical solutions
of the Riccati equation with coefficients satisfying integral or
differential conditions with arbitrary functions, Universal Journal of
Applied Mathematics \textbf{2}, 109-118 (2014).
\end{thebibliography}
\end{document}